\documentclass[12pt]{article}

\usepackage{amssymb,amsmath,amscd}

\newtheorem{df}{Definition}[section]
\newtheorem{thm}{Theorem}[section]

\newtheorem{rem}{Remark}[section]

\def\v{\vert}
\def\eq{\eqno}
\def\ot{\otimes}
\def\si{\sigma}
\def\lm{\lambda}

\begin {document}

\vspace{3cm}
\begin{center}
{\Large \bf Iwasawa decomposition of the Lie supergroup $SL(n,m,\mathbb{C})$}\\
[50pt]{\small
{\bf F. Pellegrini}\,\footnote{e-mail: pelleg@iml.univ-mrs.fr}
\\ ~~\\Institut de math\'ematiques de Luminy,
 \\ 163, Avenue de Luminy, 13288 Marseille, France}
\end{center}

\vspace{1cm}

\centerline{\textbf{Abstract}}
We show that the superanalogue of the Iwasawa decomposition exists for 
supergroup $SL(n,m,\mathbb{C})$. The first component of the decomposition is the 
compact real form $SU(n,m)$, which was constructed following the idea of 
our article [3]. The second component is the super version of 
the $AN$ group, that we define in this article.

\section{Introduction}

Iwasawa decomposition of a simple complex Lie group $D$ has many useful 
applications in group theory. For instance, it gives rise to 
Poisson-Lie structure on the compact real form $G$ of $D$, such that the 
dual Poisson-Lie group is identified with the subgroup $AN$ of 
$D=G^{\tiny{\mathbb{C}}}$. In fact, $D$ is said to be the Lu-Weinstein Drinfeld
double of $G$ and $AN$. One of our  motivation for superizing the Iwasawa 
decomposition is to define super Lu-Weinstein Drinfeld double.
The reader may be surprised that the superization of the Iwasawa 
decomposition has not yet been considered in the literature. The 
reason is simple: it was reported in [5] that real forms of 
complex simple Lie superalgebra are never compact. As the compacity is 
an important ingredient of the standard Iwasawa decomposition, 
this result seemed to imply that there is no super-version of 
the decomposition. Recently, we have shown in [3] that the 
definition of the real form given in [5] for the supercase was too 
restrictive. In fact, we argued that there is a notion of graded real 
form which is more flexible than [5] one. In particular, we 
show that working with this new concept, each supergroup of the series 
$OSp(2r,s,\mathbb{C})$ has precisely one compact graded real form. In this article
we extend the results of [3] in two directions. First, we show 
the existence of compact graded real form also for $SL(n,m,\mathbb{C})$ supergroup. 
Secondly, we argue that our new concept is very natural and powerful since 
it does allow to construct the super-analogue of the Iwasawa decomposition.\\
 
In section $2$, we give the definition of the notion of normal and 
graded real form. Next, we define the supergroup $SL(n,m,\mathbb{C})$. Finally,
we construct the compact graded real form $SU(n,m)$ of the supergroup 
$SL(n,m,\mathbb{C})$ following the ideas of our paper [3]. 
\\

In section $3$, we define the three ingredients used in the Iwasawa 
decomposition i.e. the superalgebra of complex "functions" on 
respectively the real supergroup $SL(n,m,\mathbb{C})$, $SU(n,m)$ and $s(AN)$. 
This last supergroup is the superization of the real group $AN$.
Next, we show the existence of the Iwasawa decomposition  
$SL(n,m,\mathbb{C})=SU(n,m) s(AN)$. Our method uses the superization of 
the Gram-Schmidt orthonormalisation of a family of vectors.


\section{Compact graded real form of $SL(n,m,\mathbb{C})$}

We first give the definitions of complex matrix Lie supergroup and of 
its real form. We illustrate these notions on the supergroup $SL(n,m,\mathbb{C})$.


\begin{df}
    A \underline{complex matrix Lie supergroup} $\mathcal{H}$ is a complex  
    superbialgebra generated by finite set of odd and even generators 
    subject to polynomials relations. Those relations are supposed to 
    generate a superideal of superbialgebra such that the quotient 
    can be given the structure of a Hopf superalgebra (i.e. the 
    antipode can be defined).
\end{df}    

\noindent Now we turn to the definition of normal and graded real form:


\begin{df}
    A \underline{normal real form} of a complex Lie supergroup $\mathcal{H}$ is 
    a pair $(\mathcal{H},\si)$ where $\sigma$ is an even map from 
    $\mathcal{H}$ to $\mathcal{H}$ such that: 
    $$(\si \ot \si) \Delta(x) = \Delta(\si(x)), \eq(1)$$
    $$\epsilon(\si(x)) = \overline{\epsilon(x)},  \eq(2)$$
    $$\si(\lm x + \mu y) = \bar{\lm} \si(x) + \bar{\mu}\si(y), \eq(3)$$
    $$\si(x y) = \si(x) \si(y) , \eq(4)$$
    $$S \circ \si \circ S \circ \si(x) =  x, \eq(5a)$$
    $$\si(\si(x)) =  x, \eq(6a)$$
    with $x , y \in \mathcal{H}$ and $\lm , \mu \in \mathbb{C}$.
    If the two last properties are replaced by the following:
    $$S \circ \si \circ S \circ \si(x) = (-1)^{\v x \v} x, \eq(5b)$$
    $$\si(\si(x)) = (-1)^{\v x \v} x, \eq(6b)$$
    then we have a \underline{graded real form} (cf. [3]).
\end{df}

\begin{rem}
The map $\si$ is the generalisation to the supergroup framework of the 
concept of star structure. The latter is well-known in the Hopf algebra 
literature, where the real form of a complex Hopf algebra is by definition 
the star structure. Thus, we have adapted to the supergroup context the notion 
of real form, as it is defined in the Hopf algebra setting. 
\end{rem}

Now, we turn to the definition of $SL(n,m,\mathbb{C})$. First, we recall the definition 
of the complex superbialgebra of formal power series 
$\mathbb{C}[x_{ij}¥], \; i,j=1,\ldots,n+m$. The coproduct and counit 
are defined on the generators by :
$$\Delta(x_{ij}¥)= 1 \ot x_{ij}¥+\sum_{k=1}^{n+m}x_{ij}¥\ot 1+x_{ik}¥\ot x_{kj}¥ 
    , \eq(7)$$
$$ \epsilon(x_{ij}¥)=0 .\eq(8)$$
\noindent Moreover, it is also enlightening to evaluate the coproduct of the 
elements $y_{ij}¥=\delta_{ij}¥+x_{ij}¥$. We have:
$$
\Delta(y_{ij}¥)= \sum_{k=1}^{n+m}¥y_{ik}¥\ot y_{kj}¥. \eqno(8bis)
$$
These maps are defined on all elements of $\mathbb{C}[x_{ij}¥]$ by the 
morphism property of $\Delta, \epsilon$.
The gradation of the generators $x_{ij}¥$ is $\v x_{ij }¥\v = \v i 
\v+ \v j \v$ where $\v i \v = 0, \v j\v =1$ for respectively $i=1\ldots 
n, j=n+1\ldots n+m$. We have the following standard Grassmann rules 
$x_{ij}¥x_{mn}¥=(-1)^{(\v i \v + \v j \v)(\v m \v + \v n \v)} x_{mn}¥x_{ij}¥$
for the product in $\mathbb{C}[x_{ij}¥]$.


\begin{df}
    The complex Lie supergroup $SL(n,m,\mathbb{C})$ or better \\
    $Hol(SL(n,m,\mathbb{C}))$ is the quotient of the superbialgebra $\mathbb{C}[x_{ij}¥]$ 
    by the ideal generated by the polynomial \footnote[1]{
    The superdeterminant of a supermatrix $\left(
    \begin{array}{cc}
       A & B   \\
       C & D
    \end{array}
    \right)$ is given by $sdet\left(
    \begin{array}{cc}
       A & B   \\
       C & D
    \end{array}
    \right)=\frac{det(A-BD^{-1}C)}{det(D)}$ where $det$ is the usual 
    determinant of matrices (see [1]).} $sdet(\delta+x)-1=0$. The 
    antipode is defined on the quotient by the following 
    superalgebra-antimorphism i.e $S(x_{ij}¥x_{mn}¥)=(-1)^{(\v i \v +\v j \v)
    (\v m \v + \v n \v)}S(x_{mn}¥)S(x_{ij}¥)$:
    $$ S(x_{ij}¥)= -\delta_{ij}¥+(\delta+x)_{ij}¥^{-1}, \eq(9)$$
    with $i,j=1 \ldots n+m$.
\end{df}

Here $(\delta+x)^{-1}$ means the inverse of the 
supermatrices which have for elements at the row $i$ and the column 
$j$:
$(\delta_{ij}¥+x_{ij}¥)$.
The definition of the inverse of a supermatrix 
$\left(
\begin{array}{cc}
       A & B   \\
       C & D
\end{array}
\right)$ 
reads (see [1]):
$$
\left(
\begin{array}{cc}
       A & B   \\
       C & D
\end{array}
\right)^{-1}=
\left(
\begin{array}{cc}
       (A-BD^{-1}C)^{-1} & -A^{-1}B(D-CA^{-1}B)^{-1}   \\
       -D^{-1}C(A-BD^{-1}C)^{-1} & (D-CA^{-1}B)^{-1}
\end{array}
\right), \eq(10)
$$
when $A,D$ are even invertible matrices and nothing is assumed for the 
odd matrices $B,C$.\\

\begin{rem}
    We shall use the notation $SL(n,m,\mathbb{C})$ and $Hol(SL(n,m,\mathbb{C}))$ 
    interchangeably. In particular,  we shall adopt the latter notation
    when we want to stress that in the super setting we deal with the
    holomorphic "functions" on the supergroup. \underline{We frenquently 
    write "functions" in inverted com} \underline{mas, the reason is that, 
    in fact, we are working with formal series on the Lie} 
    \underline{supergroup}. But, we think that the reader will understand 
    this abuse of notations in the sense that a lot of ideas of this article 
    are more natural in thinking about it as if we are working on functions on 
    some space.
\end{rem}

Now we follow our paper [3] and we equip $SL(n,m,\mathbb{C})$ with a 
\underline{graded} real form as follows:


\begin{thm}
The even antilinear superalgebra-morphism: 
$$\si(x_{ij})=(-1)^{(\v i \v + \v j \v) \v j \v}S(x_{ji}) \eq(11)$$ 
introduce the structure of \underline{graded} real form on 
$SL(n,m,\mathbb{C})$ in the sense of Definition 2.2. 
\end{thm}


\vskip1pc
\noindent \textbf{Proof:}\\

\noindent In fact, it is enough to prove the properties (1) to (6) 
just for the generators $x_{ij}¥$ because of the property of 
superalgebra-morphism of $\si$. Thus, the property (4) is clearly 
fullfilled. The antilinearity (property (3)) stems also from the 
definition of $\si$.\\

\noindent We develop respectively the expressions 
$\si \ot \si (\Delta(x_{ij}¥))$ and $\Delta(\si(x_{ij}¥))$ 
by using the equations (7) and (11). Thereby, we obtain respectively:
\begin{eqnarray*}
\si \ot \si (\Delta(x_{ij}¥))
    & = & 1 \ot S(x_{ji}¥) (-1)^{(\v i \v + \v j \v)\v j \v}
    +(-1)^{(\v i \v + \v j \v)\v j \v} S(x_{ji}¥) \ot 1 \\
    & &+(-1)^{(\v i \v + \v k \v)\v k \v+(\v j \v + \v k \v)\v j \v}
    S(x_{ki}¥) \ot S(x_{jk}¥),
\end{eqnarray*}
and
\begin{eqnarray*}
    \Delta(\si(x_{ij}¥)) & = & (-1)^{(\v i \v + \v j \v)\v j 
    \v}\Delta(S(x_{ji}¥))\\
    & = & (-1)^{(\v i \v + \v j \v)\v j \v}
    \tau( S \ot S (1 \ot x_{ji}¥+x_{ji}¥ \ot 1 +x_{jk}¥\ot x_{ki}¥))\\
    & = & (-1)^{(\v i \v + \v j \v)\v j \v} 1 \ot S(x_{ji}¥)+
    (-1)^{((\v i \v + \v j \v)\v j \v)}S(x_{ji}¥) \ot 1 + \\
    &   & (-1)^{\v j \v (\v j \v + \v k \v)+(\v i \v + \v k \v)\v k \v} 
    S(x_{ki}¥) \ot S(x_{jk}¥).
\end{eqnarray*}
Here $\tau$ is the flip with the following property $\tau(f \ot 
g )= (-1)^{\v f \v \v g \v}g \ot f $.
Hence, $\si \ot \si \Delta(x_{ij}¥)=\Delta(\si(x_{ij}¥))$, which 
corresponds to the property (1).\\

\noindent Now we turn to the property (2):
$$\epsilon(\si(x_{ij}¥))=\epsilon((-1)^{(\v i \v + \v j \v)\v j 
\v}S(x_{ji}¥))=(-1)^{(\v i \v + \v j \v)\v j \v} \epsilon(x_{ij}¥)= 0 =
\overline{\epsilon(x_{ij}¥)}.$$\\

\noindent Moreover, we have the property (5b) because:
\begin{eqnarray*}
    S(\si(S(\si(x_{ij}¥)))) & = & S(\si(S(S(x_{ji}¥))))(-1)^{(\v i \v + 
    \v j \v)\v j \v}\\
    & = & S(\si(x_{ji}¥)) (-1)^{(\v i \v + \v j \v)\v j \v}\\
    & = & S(S(x_{ij}¥))(-1)^{(\v i \v + \v j \v)\v j \v+(\v i \v + \v j 
    \v)\v i \v}\\
    & = & (-1)^{\v i \v + \v jÊ\v}x_{ij}¥.
\end{eqnarray*}
Here we have use the property $S(S(x))=x,\; \forall x \in SL(n,m,\mathbb{C})$.\\

\noindent Finally, it remains to prove $\si(\si(x_{ij}¥))=(-1)^{\v i \v + \v j 
\v}x_{ij}¥$.
We need to evaluate $\si(S(x_{ji}¥))$. For this we first develop the following 
identity:
$$\delta_{ij}¥ = \si((\delta_{ik}¥+x_{ik}¥)(\delta_{kj}¥+S(x_{kj}¥))) \; ,$$
thus from (11) and the antilinearity of $\si$ we deduce:
$$\delta_{ij}¥(-1)^{\v j \v} = (\delta_{kj}¥+\si(S(x_{kj}¥)))
 (\delta_{ik}¥+S(x_{ki}¥))(-1)^{\v k \v \v j \v}$$
Next, we multiply each menbers of the last equation by the right 
by $(\delta_{ip}¥+x_{ip}¥)$ and we obtain:
$$(\delta_{kj}¥+\si(S(x_{kj}¥)))\delta_{kp}¥(-1)^{\v k \v \v j \v} = 
(-1)^{\v j \v} (\delta_{jp}¥+x_{jp}¥).$$
We deduce:
$$\si(S(x_{pj}¥)) = (-1)^{(\v p \v + \v j \v) \v j \v}x_{jp}¥.$$
Thus, from this last identity we have:
$$\si(\si(x_{ij}¥))=\si(S(x_{ji}¥))(-1)^{(\v i \v + \v j \v) \v j \v}
=(-1)^{(\v i \v + \v j \v) \v j \v+(\v i \v + \v j \v) \v i 
\v}x_{ij}¥=(-1)^{\v i \v + \v j \v}x_{ij}¥,$$
which ends up the proof.
$\qquad\qquad\qquad\qquad\qquad\qquad\qquad\qquad\qquad\qquad \blacksquare$

\begin{rem}
    If $m=0$ (i.e there are no odd genererators), our definition 
    reduces to the standard compact real form $SU(n,\mathbb{C})$ of $SL(n,\mathbb{C})$.
    For this reason, we call the pair $(Hol(SL(n,m,\mathbb{C}),\si))$ the 
    compact graded real form of $SL(n,m,\mathbb{C})$. With a slight 
    ambiguity of notations, we note the compact graded real form of 
    $SL(n,m,\mathbb{C})$ by $SU(n,m)$. 
\end{rem}


\section{Iwasawa decomposition of $SL(n,m,\mathbb{C})$}

Here, we turn to the definition of the main actors which occur in 
the Iwasawa decomposition of $SL(n,m,\mathbb{C})$. We first recall the 
Iwasawa decomposition (see [2]) in the non-supercase. When $(m=0)$, 
the Iwasawa decomposition is the statement that the group $SL(n,\mathbb{C})$ 
viewed as a \underline{real group} can be decomposed as $SL(n,\mathbb{C})=GAN$
\footnote[2]{More precisely, it exists two unique maps $g,b$ from 
$SL(n,\mathbb{C})$ in (respectively) $SU(n,\mathbb{C}),AN$ such that for all $d \in 
SL(n,\mathbb{C})$ we have $d=g(d)b(d)$.} where $G=SU(n,\mathbb{C})$ is the compact real 
form of $SL(n,\mathbb{C})$ and $AN$ is the real subgroup of $SL(n,\mathbb{C})$ of upper 
triangular matrices with real positive elements on the diagonal and 
determinant equal to one. In the supercase, we need to work  with dual
objects i.e. "functions" on the supergroup. So we need: 1) the Hopf 
superalgebra $Fun(SL(n,m,\mathbb{C}))$ of "functions" on $SL(n,m,\mathbb{C})$ viewed
as \underline{a real supergroup}, 2) the Hopf superalgebra 
$Fun(SU(n,m))$ of "functions" on $SU(n,m)$ the compact graded real form of 
$SL(n,m,\mathbb{C})$ and 3) the Hopf superalgebra  $Fun(s(AN))$ of "functions" on 
$s(AN)$ which is the superization of the previous Lie group $AN$. We discuss 
the three ingredients separately.\\

\noindent 1) For finding $Fun(SL(n,m,\mathbb{C}))$ we borrow some 
inspiration from the [4] paper where the (non-super) 
q-analogue of the Iwasawa decomposition was considered. Thus the space 
of complex "functions" on the \underline{real} supergroup $SL(n,m,\mathbb{C})$ 
is $Hol(SL(n,m,\mathbb{C})) \otimes_{\mathbb{C}}¥ \overline{Hol}(SL(n,m,\mathbb{C}))$.
$\overline{Hol}(SL(n,m,\mathbb{C}))$ is a "copy" of $Hol(SL(n,m,\mathbb{C}))$ 
where the generators $x_{ij}¥$ are named $x^{\ddagger}_{ij}¥$. The 
first copy of the tensor product corresponds to "holomorphic functions" 
on $SL(n,m,\mathbb{C})$ while the second copy to "antiholomorphic functions" 
on $SL(n,m,\mathbb{C})$. Note that we use the notation $Hol(SL(n,m,\mathbb{C}))$  for
$SL(n,m,\mathbb{C})$ viewed as the \underline{complex} supergroup and the 
notation $Fun(SL$\\ $(n,m,\mathbb{C}))$ for $SL(n,m,\mathbb{C})$ viewed as the 
\underline{real} supergroup. So, we have the following definition:


\begin{df}
    The space of complex "functions" $Fun(SL(n,m,\mathbb{C}))$ on the real 
    supergroup $SL(n,m,\mathbb{C})$ is  the Hopf superalgebra:
    $$Fun(SL(n,m,\mathbb{C}))=Hol(SL(n,m,\mathbb{C})) \otimes_{\mathbb{C}}¥ 
    \overline{Hol}(SL(n,m,\mathbb{C})). \eq(12)$$  
\end{df}

\noindent 2) The construction of the Hopf superalgebra of "functions" 
on the compact graded real form of the supergroup $OSp(2r,s,\mathbb{C})$ 
was performed in detail in [3]. Here we adapt this 
construction to the compact graded real form of $SL(n,m,\mathbb{C})$, thud the 
Hopf superalgebra $Fun(SU(n,m))$ is the quotient: 
$Hol(SL(n,m,\mathbb{C})) \ot \overline{Hol}(SL(n,m,\mathbb{C}))/I$. Here $I$ is the 
superideal generated by the polynomial equations
$\si(x_{ij}¥)-x^{\ddagger}_{ij}¥=0$ ($\si$ is the map of the theorem 
2.1). The fact that $I$ is a Hopf superideal is a consequence of the 
relations (1-6) and the properties of the antipode for $Fun(SL(n,m,\mathbb{C}))$. 
Thus we have the definition:


\begin{df}
    The space of complex "functions" on the compact graded real form 
    $SU(n,m)$ of $SL(n,m,\mathbb{C})$ is the following Hopf superalgebra:
    $$Fun(SU(n,m))=Fun(SL(n,m,\mathbb{C}))/I.\eq(13)$$
    Here $I$ is the Hopf superideal generated by the polynomial 
    equations:
    $$\si(x_{ij}¥)-x^{\ddagger}_{ij}¥=0, \eq(14)$$ 
    with $\si(x_{ij}¥)=(-1)^{(\v i \v + \v j \v)\v j \v} S(x_{ji}¥).$
\end{df}

\noindent 3) Finaly, the definition of $Fun(s(AN))$ reads:


\begin{df}
    The superalgebra $Fun(s(AN))$ of complex "functions" on the real 
    supergroup $s(AN)$ is the following Hopf superalgebra:
    $$Fun(s(AN))=Fun(SL(n,m,\mathbb{C}))/J.$$
    Here, $J$ is the Hopf superideal generated by the following 
    polynomials relations:
    $$x_{ij}¥=x^{\ddagger}_{ij}¥=0 \; \mathrm{for}\; i>j , \eq(15)$$
    $$x_{ii}¥-x^{\ddagger}_{ii}¥=0. \eq(16)$$
\end{df}

\begin{rem}
    For $m=0$, the Hopf superalgebra $Fun(s(AN))$ reduces to the Hopf 
    algebra  $Fun(AN)$ of "functions" on the Lie group $AN$. If $m \not= 
    0$, the fact that $Fun(s(AN))$ is a good definition follows from 
    the fact that the super Iwasawa decomposition can be formulated 
    with it. 
\end{rem}

\begin{rem}
As $Fun(s(AN))$ and $Fun(SU(n,m))$ are both quotients of 
$Fun(SL(n,m,\mathbb{C}))$, we have the canonical projections 
$i$ and $j$ which are the superalgebra-morphisms:
$$i: Fun(SL(n,m,\mathbb{C})) \longrightarrow Fun(SU(n,m)),\eq(18)$$
$$j: Fun(SL(n,m,\mathbb{C})) \longrightarrow Fun(s(AN)).\eq(19)$$
Both $i$ and $j$ map an element $f \in Fun(SL(n,m,\mathbb{C}))$ 
into its respective cosets.\\
\end{rem}

\noindent Now we turn to the main theorem of this article, the Iwasawa 
decomposition of the real supergroup $Fun(SL(n,m,\mathbb{C})$:


\begin{thm}(\textbf{Iwasawa decomposition})\\ 
     \noindent There exist a unique pair ($\phi$,$\psi$) of 
     superalgebra-morphisms:
     $$\phi : Fun(SU(n,m)) \longrightarrow Fun(SL(n,m,\mathbb{C}))$$ 
     $$\psi: Fun(s(AN)) \longrightarrow Fun(SL(n,m,\mathbb{C}))$$ 
     such that:
     $$\phi(i(f_{(1)}¥)).\psi(j(f_{(2)}¥))= f, 
     \;\;\;  \forall f \in SL(n,m,\mathbb{C}).
     \eq(17)$$
     Here $\Delta(f)=f_{(1)}¥\ot f_{(2)}¥$ and $.$ is the standard 
     commutative multiplication in the superalgebra $Fun(SL(n,m,\mathbb{C}))$.
\end{thm}

\begin{rem}
    This theorem in the non-super case (m=0) is the dualisation of the 
    Iwasawa decomposition of the Lie group $SL(n,\mathbb{C})$ i.e. it gives the 
    Iwasawa decomposition on the space of complex "functions" on 
    $SL(n,\mathbb{C})$. When $m=0$, note that the maps $\phi, \psi$ are,
    respectively, the pullbacks of the maps $g, b$ defined in the 
    footnote 2 for the Lie group $SL(n,\mathbb{C})$.
\end{rem}


\noindent Before giving the proof of the theorem, we have to introduce 
some notations. We said that a column supervector $X$:
$$
X=
\left(
   \begin{array}{c}
       X_{1}¥\\
       \vdots \\
       X_{n}¥\\
       \chi_{1}¥\\
       \vdots\\
       \chi_{m}¥\\
   \end{array}
\right), \; \eq(20)
$$
is even (odd) when $X_{i}¥$ are even (odd) and $\chi_{m}¥$ are odd (even). 
We define the supertranspose of the supervector $X$ by:
$$
X^{st}=
\left(
   \begin{array}{cccccc}
      (-1)^{\v X \v}X_{1}¥ & \ldots & (-1)^{\v X \v}X_{n}¥ 
       & \chi_{1}¥ & \ldots & \chi_{m}¥
   \end{array}
\right). \; \eq(21)
$$
We use also the notation:
$$
X^{\ddagger}=\left(
   \begin{array}{c}
       X_{1}¥^{\ddagger}\\
       \vdots \\
       X_{n}¥^{\ddagger}\\
       \chi_{1}¥^{\ddagger}\\
       \vdots\\
       \chi_{m}¥^{\ddagger}\\
   \end{array}
\right).
\eq(22)
$$\\
\noindent We give the definition of the scalar product of two 
supervectors $X$, $Y$:
$$(X,Y)=X^{\ddagger^{st}}Y. \eq(23)$$
This scalar product have the following properties:
$$(X,Y)^{\ddagger}=(-1)^{(\v X \v+ \v Y \v)\v Y \v}(Y,X), \eq(24)$$
$$(X \lm,Y)=(-1)^{(\v X \v+1)\v \lm \v} \lm^{\ddagger} (X,Y), \eq(25)$$
where $\lm$ are possible odd or even polynoms of the generators 
$x_{ij}¥,x^{\ddagger}_{ij}¥$.
Moreover, the norm of a supervector $X$ is noted and defined by 
$\v\v X \v\v= \sqrt{(X,X)}$.\\

\noindent We define a $(m+n) \times (m+n)$ supermatrix $P$ by 
specifying either its entries or its column supervectors. In the 
first case, we note $p_{ij}¥$ the entries of the supermatrix at the 
ith row and jth column. In the second case, $P$ reads:
$$P=(P_{1}¥,\ldots,P_{n+m}¥), \eq(26)$$ 
where 
$$P_{i}¥=
\left(
   \begin{array}{c}
       p_{1\; i}¥\\
       \vdots \\
       p_{n\; i}¥\\
       p_{n+1\; i}¥\\
       \vdots \\
       p_{m+n\; i}¥\\
   \end{array}
\right). \; \eq(27)
$$
Finaly, we end this sequence of notations with two definitions:
\begin{df}
\noindent We say that a supermatrix $P$ with a unit 
superdeterminant is a $SU(n,m)$-supermatrix 
if its diagonal elements are normalized \footnote[3]{ A normalized 
formal serie is a formal serie where monomial of degree zero is
equal to one.} formal series and $P$ fullfills  
\footnote[4]{The supertranspose of a supermatrix $N=
\left(
\begin{array}{cc}
       A & B   \\
       C & D
\end{array}
\right)$ is $N^{st}=
\left(
\begin{array}{cc}
       A^{t} & C^{t}   \\
       -B^{t} & D^{t}
\end{array}
\right)$, where $A^{t},B^{t},C^{t},D^{t}$ is the usual transposition 
of the matrice $A,B,C,D$. Moreover, the entries of $N^{\ddagger}$ at 
the ith row and jth column are $n^{\ddagger}_{ij}¥$.} 
$P^{\ddagger^{st}}P=1$. 
\end{df}

\begin{df}
We say also that a supermatrix $Q$ with a unit 
superdeterminant is a $s(AN)$-supermatrix if its diagonal elements are 
real \footnote[5]{For the conjugation $\ddagger$, an element is 
real when it is equal to its conjugate.} normalized formal series and 
$Q$ is an upper triangular supermatrix. 
\end{df}


\vskip1pc
\noindent \textbf{Proof of the theorem 3.1:}\\

\noindent 
The proof of the theorem gets organized in two parts. In the first 
part, we explicitely describe the two superalgebra-morphisms 
$\phi$ and $\psi$. Next, we show that they fullfill the property (17). 
In the second part, in turn, we show the unicity of these maps.\\


\noindent Let $\tilde{\phi}, \tilde{\psi}$ be two 
superalgebra-automorphisms of  $Fun(SL(n,m,\mathbb{C}))$.
Consider supermatrices $M,\Phi,\Psi$ with elements 
$\delta_{ij}¥+x_{ij}¥$, $\tilde{\phi}(\delta_{ij}¥+x_{ij}¥)$ and 
$\tilde{\psi}(\delta_{ij}¥+x_{ij}¥)$, respectively. We know from 
definition 3.1 that $sdet(M)=1$. 
Denote $M_{l},\Phi_{l},\Psi_{l}$ the columns of $M,\Phi,\Psi$ 
, respectively. Set
$$\Phi_{l}=\frac{V_{l}}{\v \v V_{l}\v \v}, \eq(28)$$
where the supervectors $V_{l}¥$ are defined recursively by:
$$V_{1}¥=M_{1}¥,\;V_{l}¥=M_{l}-\sum_{k=1}^{l-1}
V_{k}\frac{(V_{k},M_{l})}{(V_{k},V_{k})},\; l=2,\ldots, n+m \; . \eq(29)
$$
It is easy to see that (28) imply that $\Phi$ is a 
$SU(n,m)$-supermatrix (cf. Definition 3.4), hence $\Phi$ is invertible.
Then, we also set:
$$\Psi=\Phi^{-1}M. \eq(30)$$
We deduce from (30) that $\Psi$ is a $s(AN)$-supermatrix (cf. Definition 
3.5).\\

\noindent The fact that $\Phi$ is a $SU(n,m)$-supermatrix imply that 
$\tilde{\phi}(I)=0$ ($I$ was defined in definition 3.2). 
Hence, $\tilde{\phi}$ gives rise to a superalgebra-morphism $\phi$ 
from $Fun(SU(n,m))$ to $Fun(SL(n,m,\mathbb{C}))$ by:
$$\phi(i(\delta_{ij}¥+x_{ij}¥))=\Phi_{ij}¥, \; 
\phi(i(\delta_{ij}¥+x^{\ddagger}_{ij}¥))=\Phi^{\ddagger}_{ij}¥. \eq(31)$$
The fact that $\Psi$ is a $s(AN)$-supermatrix imply that 
$\tilde{\psi}(J)=0$ ($J$ was defined in definition 3.3). 
Hence, $\tilde{\psi}$ gives rise to a superalgebra-morphism $\psi$ 
from $Fun(s(AN))$ to $Fun(SL(n,m,\mathbb{C}))$ by:
$$\psi(i(\delta_{ij}¥+x_{ij}¥))=\Psi_{ij}¥, \; 
\psi(i(\delta_{ij}¥+x^{\ddagger}_{ij}¥))=\Psi^{\ddagger}_{ij}¥. \eq(32)$$
Furthermore, the fact that $\phi,\psi,i,j$ are superalgebra-morphisms   
makes sufficient to prove (17) just for the generators $x_{ij}¥$ and
$x^{\ddagger}_{ij}¥$. Now from eq.(30) we deduce $M=\Phi \Psi$. Hence, we 
have also $M^{\ddagger}=\Phi^{\ddagger} \Psi^{\ddagger}$. 
These two last equalities give directly the validity of (17) for the 
morphisms $(\phi,\psi)$ defined by (31), (32). Thus the existence 
of $(\phi,\psi)$ is proved.\\


\noindent The reader may wish to understand better the origin of the 
formulas (28), (29). In fact, it is a superanalogue of the Gram-Schimdt 
procedure. We start with the family of column supervectors 
$M_{1}¥,\ldots,M_{n+m}¥$. In $(n+m)$-steps we construct a familly
$V_{1}¥,\ldots,V_{n+m}¥$ of orthogonal supervectors. More precisely, 
in the $k$-th step of the recursion we modify the $k$-th column in such a 
way that it becomes orthogonal to the $k-1$ previous columns. Finaly, 
we obtain (30) by the normalisation of the orthogonal family 
$V_{1}¥,\ldots,V_{n+m}¥$.\\


\noindent Now, we turn to the unicity of the maps $\tilde{\phi},
\tilde{\psi}$. We assume that there exist two distincts pairs of 
superalgebra-automorphisms of $Fun(SL(n,m,\mathbb{C}))$ 
$(\tilde{\phi}_{k}¥,\tilde{\psi}_{k}¥),\; k=1,2$ verifying 
$\tilde{\phi}_{k}¥(I)=0,\; \tilde{\psi}_{k}¥(J)=0$ and fullfilling (17). They 
give rise to two pairs of supermatrices $(\Phi_{k}¥,\Psi_{k}¥)$ 
defined for $k=1,2$ as follows:
$$\tilde{\phi}_{k}¥(\delta_{ij}¥+x_{ij}¥)=(\Phi_{k}¥)_{ij}¥, \; 
\tilde{\phi}_{k}¥(\delta_{ij}¥+x^{\ddagger}_{ij}¥)=
(\Phi_{k}¥)_{ij}^{\ddagger},\eq(33)$$
$$\tilde{\psi}_{k}¥(\delta_{ij}¥+x_{ij})=(\Psi_{k}¥)_{ij}¥, \; 
\tilde{\psi}_{k}¥(\delta_{ij}¥+x^{\ddagger}_{ij})=
(\Psi_{k}¥)_{ij}^{\ddagger}.\eq(34)$$
Firstly, from the fact that $\tilde{\phi}_{k}¥(I)=0$ 
(resp. $\tilde{\psi}_{k}¥(J)=0$) we deduce that the supermatrix 
$\Phi_{k}¥$ (resp. $\Psi_{k}¥$) are $SU(n,m)$-supermatrix 
(resp. $s(AN)$-supermatrix). Secondly, we have
the following equality between supermatrices:
$$\Phi_{1}¥\Psi_{1}¥=\Phi_{2}¥\Psi_{2}¥, \; \eq(35)$$
because $\tilde{\phi}_{k}¥, \tilde{\psi}_{k}¥$ fullfill the relation (17). 
So we obtain:
$$\Phi_{2}¥^{-1}\Phi_{1}¥=\Psi_{2}¥\Psi_{1}¥^{-1}. \; \eq(36)$$
Finally, we remark that $\Phi_{2}¥^{-1}\Phi_{1}¥$ is a 
$SU(n,m)$-supermatrix as a matricial product of 
$SU(n,m)$-supermatrices, whereas $\Psi_{2}¥\Psi_{1}¥^{-1}$ is a 
$s(AN)$-supermatrix as a matricial product of $s(AN)$-supermatrices. The 
unique supermatrix which is both a $SU(n,m)$-supermatrix and 
$s(AN)$-supermatrix is the unit supermatrix. Hence, we deduce 
that:
$$\Phi_{1}¥=\Phi_{2}¥,\; \Psi_{1}=\Psi_{2}. \; \eq(37)$$
The unicity is therefore proved.
$\qquad\qquad\qquad\qquad\qquad\qquad\qquad\qquad\qquad\quad \blacksquare$

\section{References}

\vskip1pc
\noindent [1] F.A.Berezin,
\textit{Introduction to superanalysis}, edited by A.A.Kirillov, MPAM,
Reidel Publishing
Company, Holland (1984).

\vskip1pc
\noindent [2] Anthony.W. Knapp, 
\textit{Representation theory of semisimple groups
an overview based on examples}, Princeton University Press,
Princeton mathematical series (1986).

\vskip1pc
\noindent [3] F.Pellegrini, \textit{Grassmann real form of $OSp(2r,s,\mathbb{C})$}, 
           math.RA/0311240.

\vskip1pc
\noindent [4] P. Podle\'s and S.L. Woronowicz,  
\textit{Quantum Deformation of Lorentz Group.} 
Commun. Math. Phys. {\bf 130}, 381-431 (1990).

\vskip1pc
\noindent [5] V.V.Serganova, \textit{Classification of real simple Lie 
superalgebras and symmetric superspaces}, Functional Analysis {\bf 17 n3} 
(July-September 1983) 46--54.

\end{document}